\documentclass[a4paper
	       ]{scrartcl} 
\usepackage{typearea}
\typearea[current]{calc}  

\usepackage[latin1]{inputenc}
\usepackage[ngerman, english]{babel}

\usepackage{amssymb}
\usepackage{amsmath}
\usepackage{amsthm}
\usepackage{mathtools}

\usepackage{layout}

\title{The Pólya sum kernel and Bayes estimation}
\author{Mathias Rafler\thanks{rafler@ma.tum.de}\\
}

\renewcommand{\Pr}{\mathbb P} 

\newcommand{\Poi}{\mathbf P\!\!} 

\newcommand{\Poy}{\mathsf S\!} 

\newcommand{\Pp}{\mathsf P} 
\newcommand{\GP}{\mathsf D} 

\newcommand{\B}{\mathcal{B}} 
\newcommand{\Bbd}{\mathcal{B}_{0}} 

\newcommand{\Esig}{\mathcal{E}}

\def\fsm{\text{-a.s.}}

\newcommand{\MX}{\mathcal M(X)}

\newcommand{\Mpm}{\mathcal M^{\cdot\cdot}}
\newcommand{\MpmX}{\mathcal M^{\cdot\cdot}(X)}

\newcommand{\R}{\mathbb R}

\renewcommand{\d}{\mathrm{d}}
\renewcommand{\exp}{\operatorname{exp}}
\newcommand{\e}{\operatorname{e}}

\renewcommand{\phi}{\varphi}
\renewcommand{\theta}{\vartheta}

\newcommand{\eqa}[1]{\begin{align*}#1\end{align*}}
\newcommand{\equ}[1]{\begin{equation*}#1\end{equation*}}

\newcommand{\equn}[1]{\begin{equation}#1\end{equation}}

\theoremstyle{definition}
\newtheorem{defdefinition}{Definition}[section]

\theoremstyle{plain}
\newtheorem{defsatz}[defdefinition]{Theorem}
\newtheorem{defsatzdef}[defdefinition]{Theorem and Definition}
\newtheorem{defprop}[defdefinition]{Proposition}
\newtheorem{deflemma}[defdefinition]{Lemma}
\newtheorem{deffolgerung}[defdefinition]{Corollary}

\theoremstyle{remark}
\newtheorem{defbemerkung}[defdefinition]{Remark}

\newcommand{\satz}[1]{\begin{defsatz}#1\end{defsatz}}

\newcommand{\bem}[1]{\begin{defbemerkung}#1\end{defbemerkung}}
\newcommand{\lemma}[1]{\begin{deflemma}#1\end{deflemma}}

\newcommand{\korollar}[1]{\begin{deffolgerung}#1\end{deffolgerung}}

\numberwithin{equation}{section}

\begin{document}

\maketitle

\begin{abstract}
We consider a particular Cox process from a Bayesian viewpoint and show that the Bayes estimator of the intensity measure is the so-called Pólya sum kernel, which occurred recently in the context of the construction of the so-called Papangelou processes. More precisely, if the prior, the directing measure of the Cox process, is a Poisson-Gamma random measure, then the posterior is again a Poisson-Gamma random measure and the Bayes estimator of the intensity is the Pólya sum kernel. Moreover, we extend this result to doubly stochastic Poisson-Gamma priors and give conditions under which one can identify the Bayes estimator for the intensity.
\end{abstract}

\section{Introduction}

Given some statistical model of point processes, the interest lies in deriving statements about unknown parameters from observations, which are point configurations. For Poisson processes on an Euclidean space it is possible to determine the intensity measure among the stationary ones from a single observation~\cite{NZ76}. In a Bayesian context, one starts with a probability distribution on the set of parameters, which may be interpreted as prior information or a degree of belief on the set of models, and is interested in firstly the law of the parameter given some observations, the posterior law, and secondly in the estimator for the parameter.

Staying in the context of Poisson processes, any choice of a prior distribution on the set of stationary intensity measures leads to degenerate posteriors in the sense that they are concentrated on a single intensity measure. For a larger parameter set, when such perfect estimates are not available, one is interested in finding suitable sets of priors, as discussed e.g. in~\cite{cA74}: Desirable properties are analytical tractability in the sense that it should be possible to determine the posterior law given some observation analytically, and that together with the prior, the posterior should belong to the same class of distributions. In this case the set of priors is said to be closed under sampling or conjugate.

The questions considered in~\cite{NZ76} are strongly connected to Bayesian statistics, they were discussed in an abstract form in~\cite{hF75,eD78}: In terms of point processes one starts with a consistent family of local specifications, i.e. local laws, and aims at constructing firstly all stochastic fields which are specified by this family, and secondly integral representations of these point processes as mixtures of certain extremal elements. Once the integral representation is obtained, interpreting the mixing measure as a prior, a single observation is sufficient to determine the particular extremal element. this property was called ergodic decomposability in~\cite{GW82} and is equivalent to the posterior being a Dirac measure for almost every observation.

Returning to the conjugated classes of priors, one important example is the class of Gamma distributions as priors for the family of Poisson distributions: If the prior is $\Gamma(a,r)$ distributed, then the posterior given the observation $k$ is $\Gamma(a+1,r+k)$ distributed. Hence the family of Gamma distributions is closed under sampling for Poisson models. Of course, this does not touch the question whether the family of Gamma distributions is a natural choice.

We generalize this result to general random measures and point processes: Starting with a Poisson process and its intensity measure being a Gamma-Poisson random measure, the posterior is again a Gamma-Poisson random measure. More precisely, if the prior is the Gamma-Poisson random measure with parameters $a\in\R_+$ and $\rho$ a $\sigma$-finite measure, then the posterior given the observation $\mu$, a point configuration with possibly multiple points, has parameters $a+1$ and $\rho+\mu$. In particular the class of Gamma-Poisson random measures is such a big class that there is no chance to get more precise information about the directing intensity measure from a single observation. Moreover we show that the Bayes estimator for the underlying intensity is $z(\rho+\mu)$, where $z$ depends on $a$ in a unique way.

A similar generalization from Dirichlet distributions for multinomial modells are Dirichlet processes for a point processes realizing a fixed number of points~\cite{bS96}.

In fact, the subsequently presented ideas reverse the original considerations of Zessin~\cite{hZ09}: Is there a point process which places its points according to a Pólya urn mechanism instead of  simple urn mechanism like the Poisson process. Zessin constructed this so-called Pólya sum process as the unique point process with Papangelou kernel $z(\rho+\mu)$, that is that the Pólya sum process is the unique solution of a particular integration-by-parts formula. Intuitively, the Papangelou kernel can be understood as a conditional intensity at which points are placed given the observation $\mu$. In~\cite{mR11jcma} this Pólya sum process was identified as a Cox process, the underlying random intensity being a Gamma-Poisson random measure. This connection is subject of section~\ref{sect:polyabayes}.

In a further step we allow the parameters being random and then consider Cox processes directed by doubly stochastic Gamma-Poisson random measures. Equivalently, we consider doubly stochastic Pólya sum processes. We put the results in~\cite{mR11jtp} into the Bayesian context: Point processes with local laws given by conditioned Pólya sum processes were identified as certain mixed Pólya sum processes, hence again one has an integral representation applicable for Bayesian analysis: Firstly, for any prior which is concentrated on a certain parameter set, the posterior will be concentrated on a single point and therefore determines the parameters uniquely, hence Bayes kernel and estimator of the intensity can be identified explicitly. Secondly, using the Cox representation of the Pólya sum process, we obtain the results of the first part with a doubly stochastic Poisson-Gamma process as directing measure, a result in the spirit of~\cite{cA74} for mixed Dirichlet processes. 

Implications of these results are that for particular cases we obtain that the doubly stochastic Pólya sum process is again a Papangelou process, and moreover that each point process satisfying the same integration-by-parts formula must be a doubly stochastic Pólya sum process. These ideas are presented in section~\ref{sect:dbstpolya}, the proofs are contained in section~\ref{sect:proofs}.

\section{Some random measures, point processes and results\label{sect:polyabayes}}

Let $X$ be a polish space and denote by $\B=\B(X)$ its Borel sets as well as by $\Bbd=\Bbd(X)$ the ring of bounded Borel sets of $X$. Furthermore let $\MX$ and $\MpmX$ be the space of locally finite measures and locally finite point measures on $X$, respectively, each of which is vaguely polish, the $\sigma$-algebras generated by the evaluation mappings $\zeta_{B}(\mu)=\mu(B)$, $B\in\B_{0}$. $\Mpm(X)$ is the set of observable point configurations, i.e. locally finite subsets of $X$ with possibly multiple points, and $\MX$ is the set of 'mass distributions' with finite mass in each bounded set. $\zeta_B(\mu)$ counts the number of points of $\mu\in\MpmX$ inside $B$ considering possible multiplicities.

A probability measure $\Pp$ on $\MX$ is a random measure, and if $\Pp$ is concentrated on $\MpmX$, a point process. Finally denote by $F(X)$ be the set of bounded, non-negative and measurable functions on $X$ and $F_b(X)\subset F(X)$ the subset of those functions in $F(X)$ with bounded support. Following the notation $\zeta_B(\mu)=\mu(B)$ denote by $\zeta_f(\mu)=\mu(f)=\int f\d\mu$ for $f\in F(X)$ and a measure $\mu\in\MX$ the evaluaton mapping of $f$ at $\mu$.

Apart from the finite dimensional distributions, to characterize a random measure or point process $\Pp$ uniquely, one may either use the Laplace transform
\equ{
  L_\Pp(f)=\int \e^{-\zeta_f}\d\Pp, \qquad f\in F(X)
}
or the the Campbell measure
\equ{
  C_{\Pp}(h)=\iint h(x,\mu)\mu(\d x)\Pp(\d\mu),\qquad h\in F(X\times\Mpm(X)).
}
The Campbell measure admits under suitable assumptions two disintegrations, one with respect to its intensity measure yielding Palm kernels, and one with respect to $\Pp$ itself yielding the Papangelou kernels. The famous example is the Poisson process $\Poi_\rho$ with intensity measure $\rho\in\MX$, for which
\equn{ \label{eq:rmppr:poisson}
  C_{\Poi_\rho}(h)=\iint h(x,\mu+\delta_x)\rho(\d x)\Poi_\rho(\d\mu),\qquad h\in F(X\times\Mpm(X)).
}
Moreover, there is exactly one point process which satisfies equation~\eqref{eq:rmppr:poisson}, known as Mecke's characterization of the Poisson process. The kernel $\eta(\mu,\d x)=\rho(\d x)$ does not depend on the configuration $\mu$, meaning that the point process places points independent of all other points with the same distribution.

In analogy to inductive rules defining a law given some observation, one should read equation~\eqref{eq:rmppr:poisson} as given any observed point configuration $\mu$, the intensity for another point is given by the measure $\rho$; and moreover, if some point process is given by this rule, it must be the Poisson process. In a similar way the following modification has to be understood with additional rewards to observed points. Note that neither every choice of a rule specifies a point process, nor in case of existence this point process must be unique.

Recently, Zessin considered in~\cite{hZ09} the Papangelou kernel
\equ{
  \eta(\mu,\d x)= z\bigl(\rho+\mu\bigr)(\d x)
}
for a measure $\rho\in\MX$ and some $z\in(0,1)$, which replaces the urn mechanism with replacement of the Poisson process by a Pólya urn mechanism. Instead of a point being placed according to the intensity $\rho$ independent of a present configuration, here points in the configuration $\mu$ get a reward for the intensity of the following point. Zessin answered the question of the existence and uniqueness of a point process satisfying the functional equation
\equn{ \label{eq:rmppr:polyasum}
  C_\Pp(h) = z\iint h(x,\mu+\delta_x)\bigl(\rho+\mu\bigr)(\d x)\Pp(\d\mu),\qquad h\in F.
}
Again, there is exactly one solution, the Pólya sum process $\Poy_{z,\rho}$. This point process has like the Poisson process independent increments and is infinitely divisible. In contrast to the Poisson process, $\Poy_{z,\rho}$ is not a simple point process even if $\rho$ is a diffuse mesaure. Moreover, $\Poy_{z,\rho}$ has a representation as a Cox process~\cite{mR11jcma} with its underlying random intensity measure being a Poisson-Gamma random measure, see e.g.~\cite{NZ10} for the latter process. 
More precisely, if $\GP_{z,\rho}$ is the infinitely divisible random measure with its Levy measure $\chi\bigl(\rho\otimes\tau_z\bigr)$ being the image of the product of $\rho$ and
\equ{
  \tau_z(\d r)=\frac{1}{r}\e^{-\frac{1-z}{z}r}\d r,\qquad r>0
}
under the mapping $\chi:X\times\R_+\to\MX$, $(x,r)\mapsto r\delta_x$, then the Cox representation of the Pólya sum process is
\equn{ \label{eq:polyabayes:cox}
   \Poy_{z,\rho}=\int \Poi_\kappa \GP_{z,\rho}(\d\kappa).
}
The parameter $a$ in the introduction and the parameter $z$ are linked via the relation $a=\tfrac{1-z}{z}$.


We interpret $\GP_{z,\rho}$ as a prior for the measurable family of Poisson processes. Following~\cite{GW82}, let the Bayes kernel $B$ from $\MpmX$ to $\MX$ be
\equ{
  B(\mu,H)=\frac{\d\int_H \Poi_\kappa(\,\cdot\,)\GP_{z,\rho}(\d\kappa)}{\d\Poy_{z,\rho}(\,\cdot\,)} (\mu),
}
or equivalently be defined via
\equn{ \label{eq:intro:bayeskern}
  \iint f(\mu,\kappa)\Poi_{\kappa}(\d\mu)\GP_{z,\rho}(\d\kappa)
     =\iint f(\mu,\kappa)B(\mu,\d\kappa)\Poy_{z,\rho}(\d\mu).
}
\satz{ \label{thm:rmppr:bayeskern}
Let $z\in(0,1)$ and $\rho\in\MX$. Then the posteriori measure of the Pólya sum process $\Poy_{z,\rho}$ is a Poisson-Gamma random measure, more precisely
\equ{
  B(\mu,\,\cdot\,)=\GP_{\frac{z}{1+z},\rho+\mu}.
}
}
Note that if we write $z'=\tfrac{z}{1+z}$, then there corresponds some $a'$ to $z'$ and one has $a'=a+1$. Thus the transformation of the parameters is in the same spirit as in the random variable case. Remark that $B(\mu,\,\cdot\,)$ is the superposition of two random measures
\korollar{
Under the assumptions of Theorem~\ref{thm:rmppr:bayeskern}, we have
\equ{
  B(\mu,\,\cdot\,)=\GP_{\frac{z}{1+z},\rho}\ast\GP_{\frac{z}{1+z},\mu}.
}
}
Even more, since $\mu$ is a point measure, the second measure itself is a superposition of random measures concentrated on single points.

Knowing $B(\mu,\,\cdot\,)$, the Bayes estimator of the (random) intensity measure for the Poisson process given the observation $\mu\in\MpmX$ is
\equ{
  \int \kappa B(\mu,\d\kappa),
}
we immediatly identify the Bayes estimator as the intensity measure of $B(\mu,\,\cdot\,)$.
\korollar{ \label{thm:rmppr:bayesschaetzer}
Let $z\in(0,1)$ and $\rho\in\MX$. Then the Bayes estimator for the intensity measure given an observation $\mu\in\MpmX$ is
\equ{
  b(\mu)=z\bigl(\rho+\mu\bigr)\qquad \Poy_{z,\rho}\fsm
}
}
Hence the Bayes estimator is exactly the Pólya sum kernel. The sum exactly reflects the representation of $B(\mu,\,\cdot\,)$ as the convolution of two random measures. In such a case the intensity measure necessarily is the sum of the two intensity measures. Thus if one ignored~\eqref{eq:polyabayes:cox} and denotes the point process on the rhs. by $\Pp$, then for any non-negative, measurable function $h$
\eqa{
  C_\Pp(h)&= \iiint h(x,\mu)\mu(\d x)\Poi_{\kappa}(\d\mu)\GP_{z,\rho}(\d\kappa)\\
    &= \iiint h(x,\mu+\delta_x)\kappa(\d x)\Poi_{\kappa}(\d\mu)\GP_{z,\rho}(\d\kappa)\\
    &= \iiint h(x,\mu+\delta_x)\kappa(\d x)B(\mu,\d\kappa)\Pp(\d\mu)\\
    &= \iint h(x,\mu+\delta_x)z\bigl(\rho+\mu\bigr)(\d x)\Pp(\d\mu)
}
by applying firstly Mecke's characterization of the Poisson process and~\eqref{eq:intro:bayeskern}. Since the last equation has the unique solution $\Poy_{z,\rho}$, $\Pp$ must be the Pólya sum process.
\bem{
\begin{enumerate}
  \item Instead of a fixed $z\in(0,1)$ one may start with a measurable function $z:X\to(0,1)$. As long as $z$ is bounded away from 1, calculations do not need further justifications.
  \item If $X$ is countable, $\rho$ is the counting measure on $X$ and $z$, with the just mentioned restrictions, is a probability measure on $X$, then we are in the case of the Bose gas as in~\cite{BZ11}. Given an observation of particles, the Bayes estimator of the intensity of particles occurring in certain states is the special Pólya sum kernel.
\end{enumerate}
}

\section{Doubly stochastic Pólya sum processes\label{sect:dbstpolya}}

As indicated in the introduction, the construction of H-sufficient statistics in~\cite{eD78} for a certain set of probability measures $C$ fits into the context of Bayesian statistics and has implications in connection with~\cite{mR11jcma}. At first we briefly describe the considered problem and the results obtained.

Denote by $C$ the set of all point processes $\Pr$ satisfying
\equn{ \label{eq:dbstpsp:local}
  \Pr(\,\cdot\,|\Esig_B)=\Poy_{z,\rho}(\,\cdot\,|\Esig_B)\qquad \Pr\fsm
}
for each bounded set $B\in\B$. $\Esig_B$ is a $\sigma$-algebra containing few information about the events inside the bounded set $B$ and full information about the events outside $B$, such that if $B'$ contains $B$, then $\Esig_{B'}$ is contained in $\Esig_B$. We will be more precise in a moment and only remark that $\rho$ is assumed to be a diffuse and infinite measure. Such a $\Pr$ is called stochastic field with local characteristic given by the rhs. of~\eqref{eq:dbstpsp:local}.

In~\cite{mR11jcma} there was constructed a stochastic kernel $Q$ from $\MpmX$ to $\MpmX$ satisfying
\equn{ \label{eq:dbstpsp:statQ}
  \Pr=\int Q(\mu,\,\cdot\,)\Pr(\d\mu)
}
for any $\Pr\in C$ and moreover that
\equ{
  Q\Bigl(\mu,\bigl\{\mu':Q(\mu',\,\cdot\,)=Q(\mu,\,\cdot\,)\bigr\}\Bigr)
}
holds. In fact, $Q$ is the common conditional probability of $\Pr$ conditioned on the asymptotic $\sigma$-algebra $\Esig_\infty=\bigcap_{B\in\Bbd}\Esig_B$ for every $\Pr\in C$. In the terminology of~\cite{eD78}, $Q$ is an H-sufficient statistic for $C$, in the terminology of~\cite{GW82}, $Q$ is a decomposing kernel and $C$ an ergodically decomposable simplex.

From~\cite{eD78} it follows that there exists a subset $\Delta\subseteq C$ of extreme points and a unique probability measure $V_\Pp$ on $\Delta$ such that
\equn{ \label{eq:dbstpsp:cox}
  \Pp=\int_\Delta P V_\Pp(\d P).
}
Comparing equations~\eqref{eq:dbstpsp:cox}~and~\eqref{eq:dbstpsp:statQ}, $Q(\mu,\,\cdot\,)$ is the posterior of $V_\Pp$ given the observation $\mu$.

So far this is the essence of~\cite{eD78}~and~\cite{GW82}, the in~\cite{mR11jcma} discussed choices of $(\Esig_B)_{B\in\Bbd}$ lead to the extremal points:
\begin{enumerate}
  \item if $\Esig_B$ admits the counts the points inside $B$ with multiplicity, then
    \equ{
	  \Delta=\bigl\{\Poy_{z,\rho}:z\in[0,1)\bigr\},
    }
  \item if $\Esig_B$ admits the counts the points inside $B$ without multiplicity, then
    \equ{
	  \Delta=\bigl\{\Poy_{z,w\rho}:w\in[0,+\infty)\bigr\},
    }
  \item if $\Esig_B$ admits the counts the points inside $B$ with and without multiplicity, then
    \equ{
	  \Delta=\bigl\{\Poy_{z,w\rho}:z\in(0,1),w\in(0,+\infty)\text{ or } z=w=0\bigr\}.
    }
\end{enumerate}
For any of these cases, $Q$ is identified as
\equ{
  Q(\mu,\,\cdot\,)=\Poy_{Z(\mu),W(\mu)\rho},
}
where $Z$ and $W$ are uniquely determined from the densities of points with and without multiplicity (in the first two cases one of them is fixed). Note that since we assumed $\rho$ to be an infinite measure, these densities are almost surely constant for each extremal point. We remain in the setup of the last case, i.e. start with any priori distribution $V$ on $(0,1)\times\MX$ and consider the \emph{doubly stochastic Pólya sum process} directed by $V$.

\bem{
For any priori distribution $V$ on $(0,1)\times\MX$, the doubly stochastic Pólya sum process $\Poy_V=\int \Poy_{z,\rho}V(\d z,\d\rho)$ is a Cox process $P_V$ directed by the doubly stochastic Poisson-Gamma random measure $\GP_V=\int \GP_{z,\rho}V(\d z,\d\rho)$.
}

\satz{ \label{thm:dbstpsp:bayes}
Let $V$  be a priori distribution on $(0,1)\times\MX$. If $V$ is concentrated on $(0,1)\times\{w\rho_0: w\in\R_+\}\cup\{(0,0)\}$ for some infinite and diffuse measure $\rho_0\in\MX$, then 
\begin{enumerate}
  \item $V$ is ergodically decomposable for the family $\{\Poy_{z,w\rho_0}:z\in(0,1),w\in(0,+\infty)\text{ or } z=w=0\}$ and the Bayes kernel for the parameters is
	\equ{
		\tilde{B}_V(\mu,\,\cdot\,)=\delta_{\bigl(Z(\mu),W(\mu)\rho_0\bigr)},
	}
	where $(Z(\mu),W(\mu))$ is the unique solution of the equations
    \equ{
      w\frac{z}{1-z}=U(\mu),\qquad -w\log(1-z)=V(\mu)
    }
    and $U(\mu)$ and $V(\mu)$ are the perfect estimates of the density of points with and without multiplicity, respectively. $Z$ and $W$ are $\Esig_\infty$-measurable.
  \item The Bayes kernel $B_V$ of the Cox process $P_V$ is given by
    \equ{
      B_V(\mu,\,\cdot\,)=\GP_{Z'(\mu),W(\mu)\rho_0+\mu},
    }
hence the Bayes estimator for the random intensity measure is 
    \equ{
      b_V[\mu]=Z(\mu)\bigl(W(\mu)\rho_0+\mu\bigr).
    }
\end{enumerate}
}
Thus under the assumptions of Theorem~\ref{thm:dbstpsp:bayes} we have shown that the Cox process $P_V$ is a Papangelou process:
\korollar{ \label{thm:dbstpsp:papangelou}
Let $V$ be a distribution on $(0,1)\times\{w\rho_0: w\in\R_+\}\cup\{(0,0)\}$ for some infinite and diffuse measure $\rho_0\in\MX$. Then $\Poy_V$ is a solution of the partial integration formula
\equn{ \label{eq:dbstpsp:papangelou}
  C_\Pp(h)=\iint h(x,\mu+\delta_x) Z(\mu)\bigl(W(\mu)\rho_0+\mu\bigr)(\d x)\Pp(\d\mu)
}
and $Z$ and $W$ are $\Esig_\infty$-measurable random variables. Moreover, any solution of~\eqref{eq:dbstpsp:papangelou} with $\Esig_\infty$-measurable $Z$ and $W$ and infinite and diffuse measure $\rho_0\in\MX$ is a mixed Pólya sum process.
}
\bem{
\begin{enumerate}
  \item Despite that a priori $B_V$ and $b_V$ depend on the probability measure $V$, due to the ergodic decomposability of $V$ they do not. Hence any choice of $V$ yields a Cox process $P_V$ (a mixed Pólya sum process $\Poy_V$) which solves equation~\eqref{eq:dbstpsp:papangelou} with a single Papangelou kernel. Hence we derived and discussed an integration-by-parts formula which does not have a unique solution.
  \item The fact that in Corollary~\ref{thm:dbstpsp:papangelou} $\Poy_V$ solves equation~\eqref{eq:dbstpsp:papangelou} is a straight forward application of the integration-by-parts formula for the Pólya sum process together with the ergodic decomposability of $V$. The latter property allows to estimate the parameters $z$ and $w$ from an infinite configuration of points. This does not depend on choice of mixed Pólya sum processes: In fact, if the Pólya sum processes are replaced by Papangelou processes, and if they are mixed with respect to an ergodically decomposable mixing measure, then an equation analogue to~\eqref{eq:dbstpsp:papangelou} holds for the mixture and the mixed Papangelou process is a Papangelou process itsself.
\end{enumerate}
}

\section{Proofs\label{sect:proofs}}

We identify $B(\mu,\,\cdot\,)$ by computing the Laplace transforms of both sides of equation~\eqref{eq:intro:bayeskern}. This shows Theorem~\ref{thm:rmppr:bayeskern} Before recall, e.g. from~\cite{NZ10}, that the Laplace transform of the Poisson Gamma random measure with Levy measure $\chi\bigl(\rho\otimes\tau_z\bigr)$ is
\equ{
  L_{\GP_{z,\rho}}(h)=\exp\left[ -\iint\left[ 1-\e^{-rh(x)}\right]\tau_z(\d r)\rho(\d x)\right], \qquad h\in F(X),
}
which by the relation $\log\Bigl[1+\frac{y}{a}\Bigr]=\int_0^\infty 1-\e^{-ry}\tau_z(\d r)=\int_0^\infty [1-\e^{-ry}]\tfrac{1}{r}\e^{-ar}\d r$, $a=\tfrac{1-z}{z}$, turns into
\equn{ \label{eq:bayes:lap-gp}
  L_{\GP_{z,\rho}}(h)=\exp\left[ -\int_X\log\left[ 1+\frac{zh(x)}{1-z}\right]\rho(\d x)\right], \qquad h\in F(X).
}
Secondly note that the Laplace transform of the Pólya sum process is
\equn{ \label{eq:bayes:lap-polya}
  L_{\Poy_{z,\rho}}(g)=\exp\left[ -\int_X\log\left[ 1+z\frac{1-e^{-g(x)}}{1-z}\right]\rho(\d x)\right], \qquad g\in F(X).
}

\lemma{
Let $z\in(0,1)$ and $\rho\in\MX$. 
Then for all $g,h\in F(X)$,
\equ{
  \iint \e^{-\zeta_{g\otimes h}(\mu,\kappa)}\Poi_{\kappa}(\d\mu)\GP_{z,\rho}(\d\kappa)
    = \exp\left[ -\int_X\log\left[ 1+\frac{z\bigl[1-\e^{-g}+h\bigr]}{1-z}\right]\d\rho\right].    
}
}

\begin{proof}
Note that the inner integral on the lhs. is the Laplace transform of the Poisson process with intensity measure $\kappa$, hence
\equ{
  \int \e^{-\zeta_g}\d\Poi_{\kappa}=\exp\left(-\int 1-\e^{-g}\d\kappa\right).
}
Thus what remains ist is the Laplace transform of $\GP_{z,\rho}$ at $1-\e^{-g}+h$. But this is by equation~\eqref{eq:bayes:lap-gp}
\equ{
  \exp\left[ -\int_X\log\left[ 1+\frac{z\bigl[1-\e^{-g(x)}+h(x)\bigr]}{1-z}\right]\rho(\d x)\right].\qedhere
}
\end{proof}

Next compute the rhs. of equation~\eqref{eq:intro:bayeskern}
\lemma{
Let $z\in(0,1)$ and $\rho\in\MX$. 
Then with $z'=\tfrac{z}{1+z}$ holds for all $g,h\in F(X)$
\equ{
  \iint \e^{-\zeta_{g\otimes h}(\mu,\kappa)}B(\mu,\d\kappa)\Poy_{z,\rho}(\d\mu)
    =\exp\left[ -\int_X\log\left[ 1+\frac{z\bigl[1-\e^{-g}+h\bigr]}{1-z}\right]\d\rho\right].
}
}

\begin{proof}
We check the Ansatz $B(\mu,\,\cdot\,)=\GP_{z',\rho+\mu}$. By equation~\eqref{eq:bayes:lap-gp},
\eqa{
  \int \e^{-\zeta_h(\kappa)}B(\mu,\d\kappa)
    &=\exp\left[ -\int_X\log\left[ 1+\frac{z'h(x)}{1-z'}\right]\bigl(\rho+\mu\bigr)(\d x)\right]\\
    &=\exp\left[ -\int_X\log\left[ 1+zh(x)\right]\bigl(\rho+\mu\bigr)(\d x)\right].
}
Therefore we get two exponentials, where for the integration with respect to $\Poy_{z,\rho}$ only the integral with respect to $\mu$ matters. But this just is the Laplace transform of $\Poy_{z,\rho}$ evaluated at $g+\log[1+zh]$, hence
\eqa{
  \iint &\e^{-\zeta_{g\otimes h}(\mu,\kappa)}B(\mu,\d\kappa)\Poy_{z,\rho}(\d\mu)\\
    &=\exp\left[ -\int\log\left[ 1+zh\right]\d\rho\right]\exp\left[ -\int\log\left[ 1+z\frac{1-\e^{-g-\log[1+zh]}}{1-z}\right]\d\rho\right]\\
    &=\exp\left[ -\int_X\log\left[ 1+\frac{z\bigl[1+h-\e^{-g}\bigr]}{1-z}\right]\d\rho\right].
    \qedhere
}
\end{proof}

Finally to get Corollary~\ref{thm:rmppr:bayesschaetzer}, note that for any $h\in F_b(X)$, the intensity $\nu^1_{\GP_{\frac{z}{1+z},\rho+\mu}}$ can be computed from the Campbell measure
\eqa{
  \nu^1_{\GP_{\frac{z}{1+z},\rho+\mu}}(h)&=C_{\GP_{\frac{z}{1+z},\rho+\mu}}(h\otimes 1)\\
     &=\iiint h(x)\e^{-r/z}\d r\bigl(\rho+\mu\bigr)(\d x)\GP_{\frac{z}{1+z},\rho+\mu}(\d\kappa)\\
	 &=z\int_X h(x)\bigl(\rho+\mu\bigr)(\d x)=b[\mu](h).
}

The first statement of Theorem~\ref{thm:dbstpsp:bayes} is a reformulation of the Theorem in~\cite{mR11jtp}. The second part follows from the observation that for a non-negative, measurable function $g$ by the application of Theorem~\ref{thm:rmppr:bayeskern} and the ergodic decomposability
\eqa{
  \int g(\mu,\kappa,z,w\rho)&\Poi_\kappa(\d\mu)\GP_{z,w\rho}(\d\kappa)V(\d z,\d w)\\
    &= \int g(\mu,\kappa,z,w\rho)\GP_{z',w\rho+\mu}(\d\kappa)\Poy_{z,w\rho}(\d\mu)V(\d z,\d w)\\
    &=\int g(\mu,\kappa,z,w\rho)\GP_{Z'(\mu),W(\mu)\rho+\mu}(\d\kappa)\Poy_{z,w\rho}(\d\mu)V(\d z,\d w).
}
Dropping the dependence of $g$ on the last two arguments, we get the second part of Theorem~\ref{thm:dbstpsp:bayes}.

\begin{proof}[Proof of Corollary~\ref{thm:dbstpsp:papangelou}]
Any mixed Pólya sum process $\Poy_V$ solves the partial integration formula~\eqref{eq:dbstpsp:papangelou}. Now assume that $\Pp$ is any solution of~\eqref{eq:dbstpsp:papangelou}, then the joint Laplace transform of $Z$, $W$ and $\Pp$ for $u,v,t\geq 0$, $f:X\to\R$ non-negative, bounded and measurable with bounded support is
\equ{
  L_{Z,W,\Pp}(u,v,tf)=\Pp\left(\e^{-uZ-vW-t\zeta_f}\right)=\Pp\left(\e^{-uZ-vW}\Pp\left(\e^{-t\zeta_f}|\Esig_\infty\right)\right)
}
by conditioning on $\Esig_\infty$. Denote by $\Pp_\infty$ the conditioned point process $\Pp$. Differentiation with respect to $t$ yields the Campbell measure of $\Pp_\infty$, which allows to identify this conditional measure, thus on the one hand
\equ{
  -\frac{\d}{\d t}L_{Z,W,\Pp}(u,v,tf) = \Pp\left(\e^{-uZ-vW}C_{\Pp_\infty}\left(f\otimes\e^{-t\zeta_f}\right)\right).
}
On the other hand,
\eqa{
  \begin{multlined}
  -\frac{\d}{\d t}L_{Z,W,\Pp}(u,v,tf) = C_\Pp\left(f\otimes\e^{-uZ-vW-t\zeta_f}\right)\\
     =\Pp\left(\e^{-uZ-vW}\iint f(x)\e^{-t\mu(f)-tf(x)} Z\bigl(W\rho+\mu\bigr)(\d x)\Pp_\infty(\d\mu)\right)
  \end{multlined}
}
for all $u, V\geq 0$. Exchanging integration and differentiation is justified since $f\otimes\e^{-uZ-vW-t\zeta_f}$ and $C_{\Pp_\infty}\left(f\otimes\e^{-t\zeta_f}\right)$ are integrable since $f$ is bounded and has bounded support. Thus $\Pp_\infty$ satisfies $\Pp$-a.s. the functional equation~\eqref{eq:rmppr:polyasum} and therefore is a Pólya sum process with the parameters given by $Z$ and $W\rho$. But then immediatly $\Pp$ is a mixture of these processes.
\end{proof}
\bibliographystyle{alpha-abbrv}


\end{document}